\documentclass[12pt]{article}

\usepackage{amssymb,amsthm,amscd,array}

\let\ssection=\section
\renewcommand{\section}{\setcounter{equation}{0}\ssection}

\newcommand{\bbR}{\mathbb{R}}

\newcommand{\bbC}{\mathbb{C}}

\newcommand{\bbN}{\mathbb{N}}
\newcommand{\bbZ}{\mathbb{Z}}

\newcommand{\bone}{{\bf 1}}

\newcommand{\Aff}{\mathrm{Aff}}
\newcommand{\aff}{\mathrm{aff}}

\newcommand{\cD}{{\mathcal{D}}}
\newcommand{\ce}{\mathrm{ce}}
\newcommand{\CE}{\mathrm{CE}}

\newcommand{\CO}{\mathrm{CO}}
\newcommand{\cE}{{\mathcal{E}}}
\newcommand{\Diff}{\mathrm{Diff}}

\newcommand{\rD}{\mathrm{D}}

\newcommand{\End}{\mathrm{End}}
\newcommand{\cF}{{\mathcal{F}}}

\newcommand{\gr}{{\mathrm{gr}}}
\newcommand{\rg}{\mathrm{g}}
\newcommand{\rG}{\mathrm{G}}

\newcommand{\rh}{\mathrm{h}}

\newcommand{\Id}{\mathrm{Id}}

\newcommand{\cQ}{{\mathcal{Q}}}

\newcommand{\rP}{\mathrm{P}}
\newcommand{\rR}{\mathrm{R}}
\newcommand{\rT}{\mathrm{T}}

\newcommand{\cS}{{\mathcal{S}}}

\newcommand{\SL}{\mathrm{SL}}
\newcommand{\Sl}{\mathrm{sl}}
\newcommand{\SO}{\mathrm{SO}}
\newcommand{\se}{\mathrm{e}}
\newcommand{\so}{\mathrm{o}}
\newcommand{\Sp}{\mathrm{Sp}}

\newcommand{\Supp}{\mathrm{Supp}}

\newcommand{\Vect}{\mathrm{Vect}}

\newcommand{\half}{\frac{1}{2}}
\newcommand{\fg}{\mathfrak{g}}
\newcommand{\fG}{\mathfrak{G}}
\newcommand{\fh}{\mathfrak{h}}
\newcommand{\fH}{\mathfrak{H}}

\begin{document}

\baselineskip=18pt


\def\a{\alpha}
\def\b{\beta}
\def\c{\gamma}
\def\d{\delta}
\def\g{\gamma}
\def\om{\omega}
\def\r{\rho}
\def\s{\sigma}
\def\vfi{\varphi}
\def\l{\lambda}
\def\m{\mu}
\def\implies{\Rightarrow}

\oddsidemargin .1truein
\newtheorem{thm}{Theorem}[section]
\newtheorem{lem}[thm]{Lemma}
\newtheorem{cor}[thm]{Corollary}
\newtheorem{pro}[thm]{Proposition}
\newtheorem{ex}[thm]{Example}
\newtheorem{rmk}[thm]{Remark}
\newtheorem{defi}[thm]{Definition}


\title{Projectively and conformally invariant star-products}

\author{C.~Duval\thanks{
Universit\'e de la M\'editerran\'ee and CPT-CNRS,
Luminy Case 907, F--13288 Marseille, Cedex~9, FRANCE;
mailto:duval@cpt.univ-mrs.fr
}
\and
A.M.~El Gradechi\thanks{
Facult\'e des Sciences, Universit\'e d'Artois, F--62307 Lens, FRANCE
and
CPT-CNRS, Luminy Case 907, F--13288 Marseille, Cedex~9, FRANCE;
mailto:amine@euler.univ-artois.fr
}
\and
V.~Ovsienko\thanks{
Institut Girard Desargues, Universit\'e Claude Bernard Lyon 1,
F--69622 Villeubanne, Cedex, FRANCE;
mailto:ovsienko@desargues.univ-lyon1.fr
}
}

\date{}

\maketitle

\thispagestyle{empty}

\begin{abstract}
We consider the Poisson algebra $\cS(M)$ of smooth functions on
$T^*M$ which are fiberwise polynomial. In the case where $M$ is
locally projectively (resp. conformally) flat, we seek the
star-products on $\cS(M)$ which are $\SL(n+1,\bbR)$
(resp. $\SO(p+1,q+1)$)-invariant. We prove the existence of such
star-products using the projectively (resp. conformally) equivariant
quantization, then prove their uniqueness, and study their main
properties. We finally give an explicit formula for the canonical
projectively invariant star-product.
\end{abstract}

\vskip1cm
\noindent
\textbf{Keywords:} Quantization, projective structures, conformal
structures, star-product, invariant theory.



\section{Introduction}\label{Intro}

The deformation quantization program initiated in the seventies \cite{BFFLS} was
aimed at defining an autonomous quantization method based on
Gerstenhaber's general theory of deformation of algebraic structures
\cite{Ger}. The original idea was to view quantum mechanics as a one-parameter
deformation of classical mechanics, more precisely, a one-parameter deformation of
the algebraic structures underlying classical mechanics.

If $P$ is a Poisson manifold, then
$C^\infty(P)$ is naturally equipped with two algebraic structures,
namely, the associative and commutative pointwise multiplication
and the Lie algebra defined by the Poisson bracket.
The deformed algebraic structure, describing
the quantum mechanical counterpart of $(C^\infty(P), \cdot,
\{\cdot,\cdot\})$ is $(C^\infty(P)[[\hbar]], \star)$, where the operation
$\star$, called {\it star-product},
is an associative (but non-commutative) product on $C^\infty(P)[[\hbar]]$
deforming the commutative multiplication in the direction of the
Poisson bracket. More precisely:
\begin{defi}
{\rm
Let $P$ be a Poisson manifold and $C^\infty(P)$ the space of smooth
complex-valued functions on $P$. A star-product on $P$ is an
associative algebra structure on $C^\infty(P)[[\hbar]]$, denoted
$\star$, and given by a linear map from
${C^\infty(P)\otimes{}C^\infty(P)}$ to $C^\infty(P)[[\hbar]]$,
extended by linearity to
$C^\infty(P)[[\hbar]]\otimes{}C^\infty(P)[[\hbar]]$,
such that
\begin{equation}
F\star{}G=F\cdot{}G+\frac{i\hbar}{2}\{F,G\}+
\sum_{r=2}^\infty{(i\hbar)^r{}B_r(F,G)}.
\label{starproduct}
\end{equation}
}
\end{defi}

In the mathematical literature $\hbar$ is a formal parameter,
whereas in physical applications $\hbar$ is Planck's constant.

There are usually three extra requirements for star-products~:
\begin{itemize}
\item[C1.]
the constant function $\bone$ is the unit of
$(C^\infty(P)[[\hbar]],\star)$, namely
$\bone\star{}F=F\star{}\bone=F$;
\item[C2.]
the star-product is symmetric, viz
$\overline{F\star{}G}=\overline{G}\star\overline{F}$;
\item[C3.]
\label{bidiff}
the bilinear maps $B_r$ are given by bidifferential operators.
\end{itemize}
Note that Condition C2 is sometimes called parity condition.
\medskip

The first reported star-product appeared in the work of
Gr{\oe}newold \cite{Gr}. It was derived from the Weyl-Wigner quantization on
$P=\bbR^{2n}$. It is nowadays more commonly known as the Moyal
star-product ; Moyal actually obtained the  Lie algebra bracket associated with
Gr{\oe}newold's star-product \cite{Moy}. This first star-product was later on
rediscovered by Vey~\cite{Vey}.

The general problem of existence of star-products was raised in
\cite{BFFLS}. Using cohomological techniques, De Wilde and Lecomte
\cite{DL1} proved the existence of star-products on any
symplectic manifold. A geometric proof of the same result together with
an algorithmic construction was obtained by Fedosov~\cite{Fed1,Fed2} (see
\cite{Wein} for a survey of this construction and \cite{OMY} for an
alternative approach).

More recently, Kontsevich proved an existence theorem for an
arbitrary Poisson manifold, giving explicit formul\ae\ for $P=\bbR^n$ \cite{Kon}.
An operadic and a quantum field theoretic interpretations of Kontsevich's result
were later on given respectively by Tamarkin \cite{Tam}, and Cattaneo and Felder
\cite{CF}.

The problem of the uniqueness of star-products is usually
studied modulo equivalence (see Section \ref{EquivDefSect}
for definitions and \cite{Gut,Del} for recent developments). However, extra
conditions can sometimes be imposed to single out a canonical star-product.
For instance, Gutt \cite{Gut2} proved that the Moyal star-product is the unique
$(\Sp(2n,\bbR)\ltimes\bbR^{2n})$-invariant and covariant star-product
on~$\bbR^{2n}$.  The notion of a $\fG$-invariant star-product, where $\fG$ is a
Lie group of Poisson automorphisms of $P$, was introduced in~\cite{BFFLS} (see
Section~\ref{InvCov} for definitions).  Existence of a $\fG$-invariant
star-product on a symplectic manifold was proved by Lichnerowicz
~\cite{Lic} for any compact Lie group~$\fG$ of symplectomorphisms.
More recently, Fedosov \cite{Fed3} constructed a $\fG$-invariant star-product on a
symplectic manifold endowed with a $\fG$-invariant symplectic connection.

In this article, we deal with cotangent bundles $P=T^*M$ equipped
with their canonical symplectic structure, and restrict considerations to the
Poisson algebra~$\cS(M)$ of smooth functions on $T^*M$ polynomial on
fibers. We furthermore assume $M$ to be a smooth {$n$-dimensional} manifold
endowed with either a projectively or a conformally flat structure,
i.e., $M$ admits a (locally defined) action of either $\SL(n+1,\bbR)$
or $\SO_0(p+1,q+1)$, the connected component of the pseudo-orthogonal
group with $n=p+q$.  The basic
example of a projectively (resp. conformally) flat manifold is~$\bbR{}\rP^n$
(resp.~$(S^p\times S^q)/\bbZ_2$).

Denote by $\fG$ either the projective or the conformal group.
We study, in the present article, $\fG$-invariant star-products on $T^*M$, where
the
$\fG$-action is the canonical lift of the natural action on the base.
Our first result, Theorem \ref{MainThm1}, establishes the uniqueness of a
$\fG$-invariant homogeneous star-product on $\cS(M)$. Our second
result, Theorem \ref{MainThm2}, proves the uniqueness of a $\fG$-invariant
star-product modulo $\fG$-equivalence and reparametrization.

Let us emphasize that we do not assume conditions C1, C2 and C3 a
priori.  It turns
out that C1 and C2 are automatically satisfied while C3 doesn't hold;
in fact, the maps $B_r$ in (\ref{starproduct}) are
pseudo-differential bilinear operators.
Our $\fG$-invariant star-products cannot be obtained by
Fedosov's or Kontsevich's constructions, as the latter
lead to bidifferential star-products.

\goodbreak

The existence of $\fG$-invariant star-products on~$\cS(M)$ is based on the
existence of a $\fG$-equivariant quantization map~\cite{LO,DLO}
(see also~\cite{DO2}). The latter is the unique (up to
normalization) isomorphism of $\fG$-modules, $\cQ_\l:\cS(M)\to\cD_\l(M)$, where
$\cD_\l(M)$ is the space of differential operators acting on tensor densities 
of degree~$\l$. Such a quantization map defines a $\fG$-invariant associative product on
$\cS(M)$ which turns out to be a star-product for $\l=\half$ as proved
in~\cite{B2,DLO}. The existence and uniqueness results of the present article
represent the deformation quantization counterparts of those obtained for
$\fG$-equivariant quantization. In both situations invariance properties ensure
uniqueness.

The pseudo-differential nature of the $\fG$-invariant star-products
has been revealed by Brylinski \cite{B2} and
Astashkevich and Brylinski~\cite{AB}. In the latter Reference, invariant
star-products on minimal nilpotent coadjoint orbits of semi-simple Lie
groups have been investigated. These results are closely related to
ours since these orbits are punctured cotangent bundles
$T^*M\!\setminus\!M$; nevertheless the Poisson algebras considered
in~\cite{AB} are smaller than $\cS(M)$. Moreover, our approach provides explicit
formul\ae\  in the projective case, answering a question raised in \cite{AB}.

The paper is organized as follows. In Section~2 we recall the notions of
invariant and equivalent star-products, and we give a short account on
equivariant quantization for cotangent bundles.
In Section~3, we define projective and conformal geometries and determine the ring of
projectively/conformally invariant linear operators on~$\cS(M)$. The
existence of $\fG$-invariant star-product on $T^*M$, along with a few of their
properties, are proved in Section~4. Section~5 contains our uniqueness theorems. In
Section~6, we give an autonomous derivation of the canonical projectively
invariant star-product on $\cS(\bbR{}\rP^n)$, based only on projective invariant
theory. Explicit formul{\ae} are then provided. We end this paper, with Section~7,
where we gather our conclusion, a discussion and a few perspectives.

\bigskip

\noindent\textbf{Acknowledgements:} It is a pleasure to thank Ranee Brylinski,
Simone Gutt, Pierre Lecomte and John Rawnsley for valuable help and
encouragement. This work was done while the second author was
visiting CPT as a d\'el\'egu\'e CNRS; he thanks CNRS for granting him
a d\'el\'egation and the Universit\'e d'Artois for
consenting a one year leave of absence. The second and third authors both thank the
CPT for hospitality.

\section{Invariant star-products and equivariant quantization}%
\label{Generalities}

In this section we introduce the general notions of invariance and covariance of
star-products with respect to a Hamiltonian action of a connected Lie group $\fG$.

\subsection{Invariant, covariant and strongly invariant star-products}
\label{InvCov}

First of all, let us give the precise definition of an invariant
star-product already mentioned in the Introduction.

\begin{defi}
{\rm
Given a Poisson action of a Lie group $\fG$
on a Poisson manifold $P$, a star-product $\star$ on $C^\infty(P)$ is
called $\fG$-invariant if
\begin{equation}\label{Ginvariance}
g^*(F\star{}G)=g^*F\star{}g^*G
\end{equation}
for all $F,G\in{}C^\infty(P)[[\hbar]]$ and $g\in{}\mathfrak{G}$.
}
\end{defi}

In the case where the $\fG$-action is Hamiltonian one has the following
supplementary notions.

\goodbreak

\begin{defi}\label{CovDef}
{\rm
Consider a Hamiltonian $\fG$-action on a Poisson manifold~$P$ with
associated
equivariant moment map $J:P\to\fg^*$, where $\fg^*$ is the dual of
the  Lie algebra $\fg$ of $\fG$.  A star-product on $P$ is called

a) $\fG$-\textit{covariant} if
\begin{equation}
J_X\star{}J_Y-J_Y\star{}J_X=i\hbar\,\{J_X,J_Y\}
\label{covStarProdDef}
\end{equation}

b) \textit{strongly $\fG$-invariant} if
\begin{equation}
J_X\star{}F-F\star{}J_X=i\hbar\,\{J_X,F\}
\label{strongStarProdDef}
\end{equation}
for all $F\in{}C^\infty(P)[[\hbar]]$ and $X,Y\in\fg$, where $J_X$ is the
Hamiltonian function on $P$ corresponding to~$X$.
}
\end{defi}

\begin{rmk}
{\rm
Note that a different terminology is sometimes
attached to this last notion in the literature.
What we call here strong $\fG$-invariance corresponds to the notion of preferred
observables in \cite{BFFLS,DS} and to Property $\rm I P_2$ in \cite{ACMP}. Beware
that, in the latter Reference, strong invariance means covariance and
invariance. }
\end{rmk}

\goodbreak

Let us now recall the following useful result.

\begin{pro}[\cite{ACMP}]\label{StrongImpliesInvCov}
If a star-product is strongly $\fG$-invariant, then it is both $\fG$-invariant and
$\fG$-covariant.
\end{pro}
\begin{proof}
Using the definition (\ref{strongStarProdDef}) of strong invariance, we write
\begin{eqnarray*}
i\hbar\{J_X,F\star{}G\}
&=&
J_X\star{}F\star{}G-F\star{}G\star{}J_X\\[8pt]
&=&
J_X\star{}F\star{}G
-
F\star{}J_X\star{}G+
F\star{}J_X\star{}G
-F\star{}G\star{}J_X\\[8pt]
&=&
i\hbar\left(\{J_X,F\}\star{}G+F\star{}\{J_X,G\}\right)
\end{eqnarray*}
which is nothing but the infinitesimal version of formula
(\ref{Ginvariance})
expressing the invariance property. The $\fG$-invariance of the star-product then
follows from the connectedness of $\fG$.

As for covariance, it is an immediate consequence of strong invariance.
\end{proof}

\begin{rmk}
{\rm
The converse of Proposition \ref{StrongImpliesInvCov} is proved in \cite{ACMP} under
the additional assumption of a transitive $\fG$-action.
}
\end{rmk}

\subsection{Equivalence, $\fG$-equivalence and reparametrization}
\label{EquivDefSect}

In the traditional classification of star-products one introduces a
notion of equivalence. Two star-products,
$\star$ and  $\star'$, are called equivalent if there exists a formal series
\begin{equation}
\label{EqMap}
\Phi=\Id+i\hbar\Phi_1+(i\hbar)^2\Phi_2+\cdots
\end{equation}
where $\Phi_i:C^\infty(P)\to C^\infty(P)$ are some linear operators,
such that
\begin{equation}
\label{EqUsl}
\Phi(F)\star{}\Phi(G)=\Phi(F\star'G).
\end{equation}
Usually, one also allows for formal changes of the parameter of deformation:
\begin{equation}
\label{EqRepar}
\m:i\hbar\mapsto{}i\hbar+\sum_{k=2}^\infty{}a_k(i\hbar)^k
\end{equation}
where $a_k\in\bbR$, in order to comply with Property~C2 from the Introduction.

\goodbreak

For $\fG$-invariant star-products it is natural to
consider the notion of $\fG$-equi\-valence.
\begin{defi}[\cite{Lic}]
{\rm
Two equivalent $\fG$-invariant star-products
are called $\fG$-equi\-valent if each map $\Phi_i$ in
(\ref{EqMap}) is $\fG$-equivariant.
}
\end{defi}

\noindent
The condition for two star-products to be $\fG$-equivalent is much stronger
than the usual condition of equivalence (see \cite{BBG} for recent developments).

\subsection{Equivariant quantization and the associated invariant star-product}

Equivariant quantization as developed in \cite{LO,DLO,DO1,DO2} applies
to cotangent bundles. {}From here on we restrict ourselves to $P=T^*M$
endowed with its canonical symplectic form.

Let $\cS(M)\subset C^\infty(T^*M)$ be the space of (complex-valued) functions on
$T^*M$ polynomial on fibers, and $\cD(M)$ be the space of linear differential
operators acting on $C^\infty(M)$.
The space $\cS(M)$ is the space of symbols of operators in $\cD(M)$; it has a
 natural grading
\begin{equation}
\cS(M)=\bigoplus_{k=0}^\infty\cS_k(M)
\label{grading}
\end{equation}
by the degree of homogeneity.

Let $\cF_\l(M)$ be the space of tensor densities of degree $\l\in\bbC$ on $M$,
i.e., the space of sections of the complex line bundle
$\left\vert\Lambda^n{}T^*M\right\vert^\l\otimes\bbC$. In local coordinates such
densities are of the form
\begin{equation}
f(x^1,\ldots,x^n)\,\vert{}dx^1\wedge\cdots\wedge{}dx^n\vert^\l
\label{localIdentification}
\end{equation}
 with
$f\in{}C^\infty(M)$.
Denote $\cD_\l(M)$ the space of linear differential operators on~$\cF_\l(M)$; it
has a natural filtration
$$
\cD_\l^0(M)\subset\cD_\l^1(M)\subset\cdots\subset\cD_\l^k(M)\subset\cdots
$$
such that $\cS(M)=\gr(\cD_\l(M))$.

\begin{defi}
{\rm
A {\it quantization map} is an invertible linear map
$$
\cQ_\l:\cS(M)\to\cD_\l(M)[\hbar]
$$
which preserves the principal symbol in the following
sense~:  for a homogeneous polynomial
$F\in\cS_k(M)$, the principal symbol of the differential operator~$\cQ_\l(F)$ is
 equal to
$(i\hbar)^kF$.
}
\end{defi}

There is a natural action of the group of
diffeomorphisms,
$\Diff(M)$, on $\cF_\l(M)$, denoted by $g_\l:\cF_\l(M)\to\cF_\l(M)$ for all
$g\in\Diff(M)$. We will rather use the corresponding action of the Lie algebra of
vector fields,
$\Vect(M)$, which is given by
\begin{equation}
L_X^\l{}f=X^i\frac{\partial{}f}{\partial{}x^i}+\l\,\frac{\partial{}X^i}{\partial{}x^i}f
\label{LieDerDens}
\end{equation}
for all $X=X^i\,\partial/\partial{}x^i\in\Vect(M)$, with the local
identification ${\cF_\l(M)\cong{}C^\infty(M)}$ made in~(\ref{localIdentification}).
(We will use Einstein's summation convention throughout this article.) Note that the
expression (\ref{LieDerDens}) is, indeed, independent of the choice of a
coordinate system. The canonical lift of the
$\Diff(M)$-action to $T^*M$ is automatically Hamiltonian with moment map $J$
 given by
\begin{equation}
J_X=\xi_i{}X^i\in\cS_1(M).
\label{IdentifMap}
\end{equation}

\begin{defi}
{\rm
Consider a Lie group $\fG\subset\Diff(M)$. A quantization map~$\cQ_\l$ is called
$\fG$-{\it equivariant} if
\begin{equation}
\cQ_\l(g^*F)=g_\l^{-1}\circ{}\cQ_\l(F)\circ{}g_\l^{}
\label{GIntertwining}
\end{equation}
for all $g\in{}\fG$ and $F\in\cS(M)$.
}
\end{defi}
The above formula plays a central r\^ole in the forthcoming
developments.
We will need its infinitesimal guise
\begin{equation}
L_X^\l\circ\cQ_\l(F)-\cQ_\l(F)\circ{}L_X^\l=Q_\l(L_X{}F)
\label{LieAlgEquivOfQ}
\end{equation}
for all $X\in\fg$, where $L_X$ stands for the canonical lift to $T^*M$ of the
fundamental
vector field associated with $X$.

{}From such a quantization map, we immediately obtain an associative product
given by
\begin{equation}
F\star_\l{}G=\cQ_\l^{-1}(\cQ_\l(F)\circ\cQ_\l(G)).
\label{AutomaticStar}
\end{equation}
Note that this product is not necessarily of the form (\ref{starproduct}).
 However, Condition C1 is automatically satisfied.

The following proposition is a direct consequence of the above definitions.
\begin{pro}
\label{Quant2Star}
If $\cQ_\l$ is a $\fG$-equivariant quantization map on $\cS(M)$, then the
associative product on $\cS(M)$ given by (\ref{AutomaticStar}) is
$\fG$-invariant.
\end{pro}

One wonders if there exists some extra condition sufficient to insure strong
{$\fG$-invariance} of the
$\fG$-invariant associative product (\ref{AutomaticStar}). The next proposition
 introduces
a natural geometric property of the quantization map that leads to the desired
 result.

\begin{pro}
If $\cQ_\l$ is a $\fG$-equivariant quantization map on $\cS(M)$, which furthermore
satisfies the following condition
\begin{equation}
\cQ_\l(J_X)=i\hbar\,L_X^\l
\label{ConditionExtra}
\end{equation}
for all $X\in\fg$, then the associative product on $\cS(M)$ given
by (\ref{AutomaticStar}) is
strongly $\fG$-invariant.
\label{ProMomMap}
\end{pro}
\begin{proof}
Let $X\in\fg$ and $F\in\cS(M)$, then, using successively
(\ref{AutomaticStar}), (\ref{ConditionExtra}), and (\ref{LieAlgEquivOfQ}), we get
\begin{eqnarray*}
J_X\star_\l{}F-F\star_\l{}J_X
&=&\left(\cQ_\l\right)^{-1}\left[\cQ_\l(J_X),\cQ_\l(F)\right]\\[6pt]
&=&\left(\cQ_\l\right)^{-1}\left[i\hbar\,L_X^\l,\cQ_\l(F)\right]\\[6pt]
&=&i\hbar\,L_X{}F\\[6pt]
&=&i\hbar\,\{J_X,F\}
\end{eqnarray*}
where the last equality stems from the definition of the moment map. The proof
 that
(\ref{strongStarProdDef}) holds is complete.
\end{proof}

\section{Projectively/conformally invariant operators}\label{Prelim}

We gather here
definitions and results that will be used throughout the paper. Those mainly
concern projective/conformal differential geometry. We will consider the Lie groups
$\fG=\SL(n+1,\bbR)$ and $\fG=\SO_0(p+1,q+1)$ together with their homogeneous spaces
$M=\bbR{}\rP^n$ and $M=(S^p\times S^q)/\bbZ_2$, respectively.
{}From here on, $\fG$ will stand for either of the two groups above
and~$M$ for either of the corresponding homogeneous spaces.
In the framework of Weyl's invariant theory \cite{Wey}, we will
introduce, for each geometry, $\fG$-invariant linear operators on $T^*M$
which will serve as our main tools.

\subsection{The projective and conformal symmetries}\label{ProjSub}

The real projective space of dimension $n$ is an $\SL(n+1,\bbR)$-homogeneous
space. More precisely, $\bbR{}\rP^n=\SL(n+1,\bbR)/\Aff(n,\bbR)$, where
$\Aff(n,\bbR)=\mathrm{GL}(n,\bbR)\ltimes\bbR^n$ is an affine subgroup of
$\SL(n+1,\bbR)$.

Let $x^1, x^2, \ldots, x^n$ be an affine coordinate system on $\bbR{}\rP^n$,
the fundamental vector fields associated with the $\SL(n+1,\bbR)$-action
on~$\bbR{}\rP^n$ are then given by~:
\begin{equation}
\frac{\partial}{\partial{}x^i},
\qquad
x^i\frac{\partial}{\partial{}x^j},
\qquad
x^ix^j\frac{\partial}{\partial{}x^j}\ ,
\label{ProjVF}
\end{equation}
with $i,j=1,\ldots,n$.
The vector fields (\ref{ProjVF}) correspond to translations,
linear transformations and inversions, respectively; they generate a flag of Lie
algebras
$$
\bbR^n\subset\aff(n,\bbR)\subset\Sl(n+1,\bbR).
$$

The sphere~$S^n$ with its canonical metric is a conformally flat manifold.
The same is true for $(S^p\times S^q)/\bbZ_2$ in the case of signature~$p-q$. Those
are homogeneous spaces $\SO(p+1,q+1)/\CE(p,q)$ where
$\CE(p,q)=\CO(p,q)\ltimes\bbR^n$ is the conformal Euclidean group,
$\CO(p,q)=\SO(p,q)\rtimes\bbR^*_+$, and $n=p+q$.

The fundamental vector fields associated with the $\SO_0(p+1,q+1)$-action
on $(S^p\times S^q)/\bbZ_2$ in an ``anallagmatic'' coordinate system are given
(see, e.g., \cite{DNF}) by
\begin{equation}
\frac{\partial}{\partial x^i},\qquad
x_i\frac{\partial}{\partial x^j}-x_j\frac{\partial}{\partial x^i},\qquad
x^i\frac{\partial}{\partial x^i},\qquad
x_jx^j\frac{\partial}{\partial x^i}-
2x_ix^j\frac{\partial}{\partial x^j}
\label{confGenerators}
\end{equation}
where $i,j=1,\ldots,n$ and where indices are
raised and lowered using the standard
metric~$g$ of~$(S^p\times S^q)/\bbZ_2$.
The vector fields (\ref{confGenerators}) correspond to
translations, rotations, homotheties and inversions, respectively;
they generate a flag of Lie
algebras
$$
\bbR^n\subset\se(p,q)\subset\ce(p,q)\subset\so(p+1,q+1).
$$

These two groups of transformations, $\fG$, define respectively the projective and
the conformal geometries; their Lie algebras, $\fg$, spanned by the vector fields
(\ref{ProjVF}) and~(\ref{confGenerators}) are  finite-dimensional maximal Lie
subalgebras of $\Vect(M)$, see \cite{LO,BL}.

We also introduce, for convenience, the notation $\fH\subset\fG$ for the affine Lie
subgroups $\fH=\Aff(n,\bbR)$ in the projective case, and $\fH=\CE_0(p,q)$ in the
conformal case. The corresponding Lie subalgebras will be denoted by $\fh\subset\fg$.

\subsection{Affine and Euclidean invariant operators}
\label{AffInvOpSec}

Since the group
$\Diff(M)$ of diffeomorphisms of $M$ admits a canonical lift to $T^*M$, let
us lift, accordingly, the action of $\fG$.
The search for $\fG$-invariant linear operators on $\cS(M)$ will be
dealt with in two stages.
We first consider the affine (resp. Euclidean) subgroup and determine the algebra
of $\Aff(n,\bbR)$-invariant (resp.
$(\SO_0(p,q)\ltimes\bbR^n)$-invariant) operators; in the next section we will then
enforce full $\fG$-invariance.

A classical result from invariant theory shows that the commutant of
$\Aff(n,\bbR)$ in $\End(\cS(M))$ is generated by the following two
operators
\begin{equation}
\cE=
\xi_i\frac{\partial}{\partial{}\xi_i},
\qquad
\rD=\frac{\partial}{\partial{}x^i}\frac{\partial}{\partial{}\xi_i}
\label{ED}
\end{equation}
which span the Lie algebra $\aff(1,\bbR)$. Indeed, an $\Aff(n,\bbR)$-invariant
linear operator mapping $\cS_k(M)$ into $\cS_\ell(M)$ is proportional to
$\rD^{k-\ell}$ (see, e.g.,
\cite{Wey,GW}). The commutant of $\Aff(n,\bbR)$ in $\End(\cS(M))$ is, hence, given
by series in $\cE$ and $\rD$, convergent on $\cS(M)$.

It has been shown in \cite{DLO} that the commutant of $\SO_0(p,q)\ltimes\bbR^n$ in
$\End(\cS(M))$ is generated by the operators
\begin{equation}
\rR=\xi^i\xi_i,\qquad
\cE=\xi_i\frac{\partial}{\partial \xi_i},\qquad
\rT=\frac{\partial}{\partial \xi^i}\frac{\partial}{\partial \xi_i}
\label{sl2}
\end{equation}
whose commutation relations are those of $\Sl(2,\bbR)$, together with
\begin{equation}
\rG=\xi^i\frac{\partial}{\partial x^i},\qquad
\rD=\frac{\partial}{\partial{}x^i}\frac{\partial}{\partial{}\xi_i},\qquad
\Delta=\frac{\partial}{\partial x^i}\frac{\partial}{\partial x_i}
\label{rh1}
\end{equation}
which span the Heisenberg Lie algebra $\rh_1$. The operators (\ref{sl2}) and
(\ref{rh1}) span the Lie algebra  $\Sl(2,\bbR)\ltimes\rh_1$.

\subsection{Projectively and conformally invariant 
operators}\label{GInvOp}

It is noteworthy that $\cE$ commutes with the lift of any diffeomorphism of $M$.
One may ask if, upon
restriction to $\fG\subset\Diff(M)$, there exist other linear operators on~$T^*M$ that
commute with $\fG$. The answer is negative in the projective case and positive in the
conformal case.

\begin{pro}\label{Projinvariants}
The commutant of $\SL(n+1,\bbR)$ in $\End(\cS(M))$ is generated by~$\cE$.
\end{pro}
\begin{proof}
An affinely invariant linear operator is a series in $\cE$ and $\rD$ of the form
\begin{equation}
A=\sum_{s=0}^\infty{}P_s(\cE)\,\rD^s,
\label{AffInvOp}
\end{equation}
where $P_s$ is a series in one variable.
Let $X_i=x^ix^j\,{\partial}/{\partial{}x^j}$ be the $i$-th generator of
inversions in (\ref{ProjVF}). Straightforward computation (see \cite{LO}) yields
the commutation relation
$$
\left[
L_{X_i},\rD
\right]
=
(2\cE+n+1)\circ\frac{\partial}{\partial\xi_i}.
$$
One then checks that
\begin{equation}
\left[L_{X_i},A
\right]
=
\sum_{s=0}^\infty\,sP_s(\cE)(2\cE+n+s)\,\rD^{s-1}\circ\frac{\partial}{\partial\xi_i}.
\label{LieDerAffInvOp}
\end{equation}
This expression vanishes if and only if $P_s=0$ for all $s\geq1$. Hence
$A=P_0(\cE)$ is a necessary condition for $A$ to commute with the
$\SL(n+1,\bbR)$-action.
\end{proof}

The conformal counterpart of the above statement is as follows.

\begin{pro}\label{Confinvariants}
The commutant of $\SO_0(p+1,q+1)$ in $\End(\cS(M))$ is the commutative
associative algebra generated by $\cE$ and the operator $\rR_0=\rR\circ\rT$.
\end{pro}
\noindent
\begin{proof}
A sketch of this proof was given in~\cite{DLO}; for the sake of
completeness we give here the details.

Let us consider an operator $Z$ on the
space of polynomials of degree $k$,
namely
$$
\cS^k(M)=\bigoplus_{\ell=0}^k\cS_\ell(M),
$$  
and commuting with the canonical lift of
$\SO_0(p+1,q+1)$. It is, according to classical invariant theory~\cite{Wey,GW}, a
differential operator, polynomial in the generators (\ref{sl2}) and (\ref{rh1}).

We therefore seek a differential operator $Z$ on $T^*M$ which commutes with
the $\SO_0(p+1,q+1)$-action. Its principal symbol $\s(Z)$ is
a function on $T^*(T^*M)$, polynomial on fibers. More precisely,
if $(\zeta_i,y^i)$ denote the conjugate variables to $(x^i,\xi_i)$
respectively, then $\s(Z)$ is polynomial in the variables $\xi_i,\zeta_i,y^i$.
The function $\s(Z)$ has to be annihilated by the canonical lifts
to $T^*(T^*M)$ of all generators (\ref{confGenerators}) of the
conformal Lie algebra $\so(p+1,q+1)$.

Let us assume that $\s(Z)$ is $\ce(p,q)$-invariant and consider then invariance 
with respect to
inversions whose $i$-th generator is
$X_i=x_jx^j\,\partial/\partial x^i-
2x_ix^j\,\partial/\partial x^j$. Its canonical lift to $T^*(T^*M)$ is given by
\begin{equation}
\label{DoubleLift}
\begin{array}{rcl}
\widetilde{L}_{X_i} & =&
\displaystyle
x_jx^j\,\frac{\partial}{\partial x^i}-
2x_ix^j\,\frac{\partial}{\partial x^j}\\[10pt]
&&
\displaystyle
+2x_i\left(
\xi_j\,\frac{\partial}{\partial\xi_j}
- y^j\,\frac{\partial}{\partial y^j}
+\zeta_j\,\frac{\partial}{\partial\zeta_j}
\right)\\[10pt]
&&
\displaystyle
-2x^j\left(
\xi_i\,\frac{\partial}{\partial\xi^j}-\xi_j\,\frac{\partial}{\partial\xi^i}
+ y_i\,\frac{\partial}{\partial y^j}- y_j\,\frac{\partial}{\partial y_i}
+\zeta_i\,\frac{\partial}{\partial\zeta^j}-\zeta_j\,\frac{\partial}{\partial\zeta^i}
\right)\\[10pt]
&&
\displaystyle
+2\left(
\xi_i\,y_j\,\frac{\partial}{\partial\zeta_j}
-y_i\,\xi_j\,\frac{\partial}{\partial\zeta_j}
-\xi_j\,y^j\,\frac{\partial}{\partial\zeta^i}
\right)
\end{array}
\end{equation}
and the invariance with respect to inversions reads $\widetilde{L}_{X_i}\s(Z)=0$.
Now, invariance with respect to $\ce(p,q)$ clearly implies that $\s(Z)$ is
annihilated by the first three terms in (\ref{DoubleLift}), so that
\begin{equation}
\label{FinalCond}
\left(\xi_i\,y_j\,\frac{\partial}{\partial\zeta_j}
-y_i\,\xi_j\,\frac{\partial}{\partial\zeta_j}
-\xi_j\,y^j\,\frac{\partial}{\partial\zeta^i}
\right)\s(Z)
=0
\end{equation}
for all $i=1,\ldots,n$.

\goodbreak

\begin{lem}
The equation (\ref{FinalCond}) implies
\begin{equation}
\label{VeryFinalCond}
\frac{\partial\s(Z)}{\!\!\!\!\!\!\partial\zeta_i}
=0
\end{equation}
 for all $i=1,\ldots,n$.
\end{lem}
\begin{proof}
The determinant of the matrix
$$
A_j^i=y^i\xi_j-\xi^iy_j+\xi_ky^k\,\d^i_j
$$
intervening in (\ref{FinalCond}) is $\det(A)=\xi_i\xi^i\,y_jy^j\,(\xi_ky^k)^{n-2}$ which is non-zero on
the complement of a lower-dimensional smooth submanifold of $T^*(T^*M)$.
\end{proof}

By $\se(p,q)$-invariance, the operator $Z$ is a polynomial in the
differential operators (\ref{sl2}) and (\ref{rh1}), see Section \ref{AffInvOpSec}.
Furthermore, invariance with respect to the generator of homotheties
$X_0=x^i\,\partial/\partial x^i$
shows that $Z$ is in fact a polynomial in 
\begin{equation}
\label{HomogOpers}
\rR_0=\rR\circ\rT,
\qquad\cE,
\qquad\rG_0=\rG\circ\rT,
\qquad\rD,
\qquad\Delta_0=\Delta\circ\rT.
\end{equation}
The principal symbols of the last three operators are
$$
\s(\rG_0)=\xi_i\zeta^i\,y_jy^j,
\qquad
\s(\rD)=\zeta_i\,y^i,
\qquad
\s(\Delta_0)=\zeta_i\zeta^i\,y_jy^j.
$$
These three polynomials are algebraically independent for $n>1$. Condition
(\ref{VeryFinalCond}) implies then that $Z$ depends only on $\cE$
and
$\rR_0$. Note that if $n=1$, we find $\rR_0=\cE(\cE-1)$ in agreement with
Proposition
\ref{Projinvariants}.

We have thus shown that, for all $k$, any $Z\in\End(\cS^k(M))$
commuting with the $\SO_0(p+1,q+1)$-action is polynomial in $\cE$ and $\rR_0$. 
This completes the proof of Proposition \ref{Confinvariants}.
\end{proof}

\section{Existence of projectively and conformally invariant
star-products}\label{ExistenceSection}

Taking advantage of the results obtained in equivariant
quantization (see \cite{LO,DLO,B2}) and of Proposition \ref{Quant2Star}, one defines a
$\fG$-invariant star-product on $T^*M$ . In this section we give a brief
account on the projectively and conformally equivariant quantizations and discuss
the main properties of the associated invariant star-products.

\subsection{Construction of $\fG$-invariant star-products}\label{ExSubsec}

It has been proved in \cite{LO,DLO} that, for any $\l\in\bbC$, there exists a
unique $\fG$-equivariant quantization map $\cQ_\l:\cS(M)\to\cD_\l(M)[\hbar]$ on
$T^*M$.

In a local coordinate system, one can locally identify $\cS(M)$ and $\cD_\l(M)$
through the \textit{normal ordering} prescription:
\begin{equation}
P^{i_1\ldots{}i_k}\xi_{i_1}\cdots\xi_{i_k}\mapsto{}
(i\hbar)^k
P^{i_1\ldots{}i_k}\frac{\partial}{\partial x^{i_1}}
\cdots\frac{\partial}{\partial x^{i_k}}
\label{NormalOrdering}
\end{equation}
where $P^{i_1\ldots{}i_k}$ is a smooth function of $(x^1,\ldots,x^n)$.

The explicit formula of $\cQ_\l$ is only known in the projective case; it is given,
in an adapted coordinate system, and using the identification
(\ref{NormalOrdering}), by the series~\cite{DO2}
\begin{equation}
\cQ_\l
=
\sum_{r=0}^\infty{
C_r(\cE)\,(i\hbar\rD)^r
}
\label{QPEbis}
\end{equation}
where $\cE$ and $\rD$ are as in (\ref{ED}) and
\begin{equation}
C_r(\cE)
=
\frac{1}{r!}
\frac{(\cE+(n+1)\l)_r}{\left(2\cE+n+r\right)_r},
\label{Coeffm}
\end{equation}
where $(a)_r:=a(a+1)\cdots(a+r-1)$ is the Pochhammer symbol. The
expression~(\ref{QPEbis}) is well defined globally on
$T^*M$ since
$M$ is projectively flat.

An important feature of the quantization map (\ref{QPEbis}) is that it is
homogeneous in the following sense. Let us assign a degree to the
 deformation
parameter $\hbar$, more precisely, we put
$$
\deg\hbar=1.
$$
Then $\cQ_\l$ preserves the total degree on
$\cS(M)[\hbar]$. In other words, one has

\begin{pro}
The quantization map
$\cQ_\l$ commutes with the Euler operator:
\begin{equation}
\widehat{\cE}=\cE+\hbar\frac{\partial}{\partial\hbar}.
\label{EulerNew}
\end{equation}
\end{pro}
\begin{proof}
This follows from the commutation relation $[\cE,\rD]=-\rD$ and the
expression~(\ref{QPEbis}).
\end{proof}

In the conformal case we have no explicit formula for the
$\SO_0(p+1,q+1)$-equivariant
quantization map. However, one can guarantee \cite{DLO} that $\cQ_\l$ is
also homogeneous in this case.

A $\fG$-invariant
star-product on $T^*M$ can be obtained from such a $\fG$-equivariant
quantization map.

\begin{pro}[\cite{DLO,B2}]
The associative product associated with $\cQ_\l$ through (\ref{AutomaticStar})
 is a star-product if and only if $\l=\half$.
\label{ProStarQ}
\end{pro}

\noindent
The proof consists in checking that $\l=\half$ is the only value of $\l$ for which
the first-order term in $\hbar$ of the
associative product~(\ref{AutomaticStar}) coincides with the Poisson bracket.

Note, however, that the uniqueness of~$\cQ_\half$ does not a priori insure
the uniqueness of a $\fG$-invariant star-product.

\subsection{Main properties}

For the constructed $\fG$-invariant star-products Condition C1 from Section
\ref{Intro} is satisfied. We will show below that Condition C2 also holds.

\begin{defi}\label{defHomogeneousStar}
{\rm
A star-product on the space $\cS(M)$ will be called \textit{homogeneous}, if all
the bilinear operators $B_r$ in (\ref{starproduct}) are homogeneous of degree $r$,
that is, if they preserve the grading (\ref{grading}) according to
\begin{equation}
B_r:\cS_k(M)\otimes\cS_\ell(M)\to\cS_{k+\ell-r}(M),
\label{HomogFor}
\end{equation}
or, equivalently, if $\widehat{\cE}$ is a derivation of the star-product
algebra.
}
\end{defi}

\goodbreak

\begin{pro}\label{SymmStarQhalf}
The $\fG$-invariant star-product (\ref{AutomaticStar}) obtained from the
$\fG$-equi\-variant quantization map $\cQ_\half$ is symmetric, homogeneous and
strongly $\fG$-invariant.
\end{pro}
\begin{proof}
The quantization map $\cQ_\half$ is symmetric, namely, it satisfies
$$
\cQ_\half(F)^*=\cQ_\half(\overline{F})
$$
for all $F\in\cS(M)$ \cite{LO,DLO,DO1}, where $\cQ_\half(F)^*$ denotes the
formal adjoint operator with respect to the natural pairing on compactly
supported $\half$-densities. Using the definition (\ref{AutomaticStar}) of the
star-product, one now gets the symmetry condition~C2.

Homogeneity of the
quantization map $\cQ_\half$ readily implies the homogeneity of the corresponding
star-product.

The projectively and the conformally equivariant quantization maps $\cQ_\half$
coincide up to second-order terms, namely, in both cases one has
$$
\cQ_\half=\Id+\frac{i\hbar}{2}\,\rD+O(\hbar^2)
$$
in any coordinate system (cf. \cite{LO,DLO}). One easily verifies that
 $\cQ_\half$
satisfies condition~(\ref{ConditionExtra}). By Proposition
\ref{ProMomMap}, the associated
$\fG$-invariant star-products are thus strongly $\fG$-invariant.
\end{proof}

Condition C3 fails to be satisfied (as proved in \cite{B2} and \cite{AB} for a
subalgebra  of~$\cS(M)$). Each term
$B_r$ is a pseudo-differential bilinear operator, while its restriction
$B_r\vert_{\cS_k(M)\otimes\cS_\ell(M)}$ is a bidifferential operator, 
just like $\cQ_\half\vert_{\cS_k(M)}$ is a differential operator, 
see~\cite{LO}.
Hence the constructed star-product is local, namely, for all
$F,G\in\cS(M)$, $\Supp(F\star{}G)\subset\Supp(F)\cap\Supp(G)$, 
see Lemma \ref{LocalityLemma} below.

\section{Uniqueness of $\fG$-invariant star-product}

Our goal is to show that the star-products
constructed in Section \ref{ExSubsec} with the help of the $\fG$-equivariant
 quantization
map are the unique $\fG$-invariant star-products where, as above, $\fG=\SL(n+1,\bbR)$
and $\fG=\SO_0(p+1,q+1)$, respectively. We prove uniqueness in two
different settings:
\begin{enumerate}
\item
in the class of homogeneous $\fG$-invariant star-products,
\item
in the class of all $\fG$-invariant star-products modulo formal reparametrizations
and $\fG$-equivalence.
\end{enumerate}

\subsection{Homogeneous star-products}\label{HomogUniq}

Homogeneity of a star-product (see Definition \ref{defHomogeneousStar}) is a very
natural property from a physical standpoint. Indeed, if one considers $\hbar$ as a
physical constant whose dimension is that of an action (i.e., the dimension of
Planck's constant which is also the inverse dimension of the Poisson bracket
on $T^*M$) then the physical dimension
of the star-product $F\star{}G$ of two observables is the same as that
of their product $FG$, when $\star$ is homogeneous. This is a direct consequence of the fact that
$B_r$ has the same physical dimension as $\hbar^{-r}$, which follows from
associativity.  In other words a homogeneous star-product is dimensionless.

On the other hand, homogeneous star-products 
were thoroughly studied in the mathematical literature. For instance, De Wilde 
and Lecomte proved \cite{DL0} that any two homogeneous star-products on a 
cotangent bundle  are equivalent (in the sense of the definitions of 
Section~\ref{EquivDefSect}). The $\fG$-invariant star-products
constructed in Section~\ref{ExSubsec} are also homogeneous (see Proposition
\ref{SymmStarQhalf}).

The first main result of this paper is
\begin{thm}
\label{MainThm1}
There exists a unique homogeneous $\fG$-invariant star-product
on the space of symbols $\cS(M)$.
\end{thm}

\begin{proof}
In Section \ref{ExistenceSection} we proved the existence of a homogeneous
$\fG$-invariant star-product on $\cS(M)$. We will now prove its uniqueness.

Let $\star$ and $\star'$ be two homogeneous $\fG$-invariant
star-products. Let us assume that the first~$r-1$ terms of these star-products
coincide, and use induction over~$r$.
The difference $B_r-B'_r$ is a $\fG$-invariant homogeneous Hochschild
$2$-cocycle.
Indeed, associativity of the star-product $\star$ implies that $\d{B_r}$
depends only upon $B_i$ with~$i<r$, where the Hochschild coboundary of a $2$-cochain
$B$ is given by
\begin{equation}
\d{B}(F,G,H)
=
F B(G,H)-B(FG,H)+B(F,GH)-B(F,G)H,
\label{coboundary}
\end{equation}
implying that $\d(B_r-B'_r)=0$.

Let $C$ be a Hochschild $2$-cocycle on $\cS(M)$. Assume now that $C$ is
homogeneous as in (\ref{HomogFor}) and $\fG$-invariant. As a bilinear map, $C$
decomposes as a sum
$C_1+C_0$, where $C_1$ and $C_0$ are, respectively, the skew-symmetric and the symmetric
parts of~$C$. We will need the following well-known result.
\begin{pro}\label{LemHoch}
For any local Hochschild $2$-cocycle $C$ on $\cS(M)$, the skew-sym\-metric part
$C_1$ is a bivector, and the symmetric part
$C_0$ is the coboundary of a local $1$-cochain.
\end{pro}
\noindent
This statement is an important result in deformation theory. It was first
established in the differentiable case \cite{Vey} and was later on generalized 
to local cocycles in \cite{CGW} (let us mention that this result also holds for
continuous cocycles \cite{Con,Nad}).

In order to apply Proposition \ref{LemHoch}, we will prove that each term $B_r$
of a $\fG$-invariant star-product is local, a result that generalizes Theorem~5.1
in~\cite{LO}.

\begin{lem}\label{LocalityLemma}
Any linear $\fG$-invariant operator $B:\cS_k(M)\otimes\cS_\ell(M)\to\cS_m(M)$ with
$m\leq{}k+\ell$ is local.
\end{lem}
\begin{proof}
We must prove that $\Supp(B(F,G))\subset\Supp(F)\cap\Supp(G)$ for all
$F\in\cS_k(M)$ and $G\in\cS_\ell(M)$. Suppose that one of the arguments, $F$ or
$G$, vanishes in a neighbourhood of some $x\in{}M$; we will prove that
$B(F,G)(x)=0$. Let us now locally identify $M$ with $\bbR^n$ and consider the
subalgebra $\bbR\ltimes\bbR^n$ of $\fg$ generated by the Euler vector field,
$\cE$, and the translations. Using translation-invariance, we may, hence, assume
$x=0$.

We will embed $\cS_k(\bbR^n)\otimes\cS_\ell(\bbR^n)$ into
$\cS_{k+\ell}(\bbR^{2n})$ and notice that $F\otimes{}G$ vanishes in a
neighbourhood of $x=0$ in $\bbR^{2n}$. It remains to show that if
$B:\cS_{k+\ell}(\bbR^{2n})\to\cS_m(\bbR^n)$ is a linear map which commutes with
the action of homo\-theties~$L_\cE$, then for all $H\in\cS_{k+\ell}(\bbR^{2n})$
that vanishes in a neighbourhood of $x=0$, we have $B(H)(0)=0$ provided
$m\leq{}k+\ell$. But the proof of the latter statement coincides with that of
Theorem~5.1 in~\cite{LO}.
\end{proof}

The building blocks of the operators $B_r$ are the $\fH$-invariant operators listed
in~(\ref{HomogOpers}). These operators never increase the degree of homogeneity in
$\xi=(\xi_i)$, hence Lemma \ref{LocalityLemma} applies. We are now able to use
Proposition~\ref{LemHoch} and consider~$C_1$ and $C_0$ separately. The assertion
of Theorem \ref{MainThm1} will follow from Lemmas
\ref{SkewLemHom} and~\ref{SymLem} below.

\begin{lem}\label{SkewLemHom}
There is no non-zero $\fG$-invariant bivector on $T^*M$ with coefficients in
$\cS(M)$ homogeneous of degree $r\geq2$.
\end{lem}
\begin{proof}
There is clearly no non-zero such bivector
$W:\cS_k(M)\otimes\cS_\ell(M)\to\cS_{k+\ell-r}(M)$, for $r>2$. For~$r=2$, if it
exists it is necessarily of the form
$
W=W_{ij}\,\partial/\partial{\xi_i}\wedge\partial/\partial{\xi_j}
$
with coefficients $W_{ij}$ of degree $0$ in $\xi$. Since $W$ is
invariant with respect to the generators of translations,
$\partial{W_{ij}}/\partial{x^s}=0$ for all
$s=1,\ldots,n$. But, in this case, $W$ cannot be invariant with respect to homotheties.

We thus have proved that there is no non-zero bivector invariant with respect to the
$(n+1)$-dimensional Lie algebra of translations and homotheties. This Lie algebra is a
Lie subalgebra of both $\Sl(n+1,\bbR)$ and
$\so(p+1,q+1)$.
Lemma \ref{SkewLemHom} is proved.
\end{proof}

\begin{rmk}
{\rm
Note that, in the proofs of Lemmas \ref{LocalityLemma} and \ref{SkewLemHom}, 
we only needed invariance with respect to a subalgebra of~$\fg$. }
\end{rmk}

\begin{lem}
There is no non-zero $\fG$-invariant Hochschild
$2$-co\-boundary $C_0$ on the
associative commutative algebra $\cS(M)$ homogeneous of degree $r\geq2$.
\label{SymLem}
\end{lem}
\noindent
\begin{proof}
Suppose that such a $C_0$ exists. Being a coboundary, it is of the form $C_0=\d{A}$ where
\begin{equation}
\d{A}(F,G)=FA(G)-A(FG)+A(F)G
\label{twoCoboundary}
\end{equation}
for some linear map  $A:\cS_k(M)\to\cS_{k-r}(M)$, with $r\geq2$.
Let us prove that $A$ is $\fG$-invariant.

Since $C_0$ is $\fG$-invariant, then, for any $X\in\fg$, the linear
map $L_X(A)=\left[L_X,A\right]$ is a Hochschild $1$-cocycle on
$\cS(M)$. Indeed one has $\d\circ{}L_X=L_X\circ\d$. Thus, $L_X(A)$ is
a derivation on $\cS(M)$. Therefore, this is a vector field on
$T^*M$ polynomial in $\xi$ and, hence, $L_X(A)$ cannot decrease the
degree by more than $1$. Note, however, that
$L_X(A):\cS_k(M)\to\cS_{k-r}(M)$ with $r\geq2$ since, again,
$L_X(A)=L_X\circ{}A-A\circ{}L_X$ and~$L_X$ preserves
$\cS_k(M)$ for any vector field $X$ on $M$. It follows that
$L_X(A)=0$ for all $X\in\fg$ and thus $A$ is $\fG$-invariant.

The classification of $\fG$-invariant linear maps on $\cS(M)$ is given by
Proposition~\ref{Projinvariants} and Proposition \ref{Confinvariants}.
Being homogeneous of degree zero in $\xi$, a non-zero
$\fG$-invariant element $A$ of
$\End(\cS(M))$ cannot decrease the degree. Lemma~\ref{SymLem}
is proved.
\end{proof}

Lemmas \ref{SkewLemHom} and \ref{SymLem} imply that $B_r-B'_r=0$ for $r\geq2$.
This completes the proof of Theorem~\ref{MainThm1}.
\end{proof}

The unique homogeneous $\fG$-invariant star-product will be called
$\fG$-\textit{canonical}. According to Proposition \ref{SymmStarQhalf}, this
$\fG$-canonical star-product is the one associated with the $\fG$-equivariant
quantization map $\cQ_\half$ from Section~\ref{ExistenceSection}. The
same Proposition also states that it is both symmetric and strongly
$\fG$-invariant.

\subsection{Uniqueness up to $\fG$-equivalence and
reparametrization}\label{UniqueUpGequiReparam}

The following theorem is the second main result of this paper.

\begin{thm}
\label{MainThm2}
The $\fG$-canonical star-product
on the space of symbols~$\cS(M)$ is the unique $\fG$-invariant star-product modulo formal
reparametrizations and $\fG$-equivalence.
\end{thm}

\begin{proof}
Let $\star$ and $\star'$ be two $\fG$-invariant star-products. Let us assume that
there exists a $\fG$-invariant formal series  (\ref{EqMap}) and a reparametrization
(\ref{EqRepar}) intertwining the first~$r-1$ terms of these star-products, and use
induction over $r$. Using this equivalence we can assume that
$\star$ and $\star'$ coincide up to the $(r-1)$-th order term. The difference $B_r-B'_r$
is then a $\fG$-invariant Hochschild $2$-cocycle.

As in Section \ref{HomogUniq} we consider the decomposition
$C=C_1+C_0$, where $C_1$ and~$C_0$ are, respectively, the skew-symmetric and the symmetric
parts of $C$. By Proposition~\ref{LemHoch}, $C_1$ is a bivector and
$C_0$ is a coboundary.

We will need the following two lemmas.
\begin{lem}\label{SkewLem}
(i) In the projective case, the canonical Poisson bivector
\begin{equation}
\Pi =
\frac{\partial}{\partial\xi_i}\wedge\frac{\partial}{\partial{}x^i}
\label{canonicalBivector}
\end{equation}
on $T^*M$ is the unique (up to an overall multiplicative constant) $\fG$-invariant bivector.

(ii) In the conformal case with $n\neq2$, the canonical Poisson bivector on $T^*M$
is the unique (up to an overall multiplicative constant) $\fG$-invariant bivector.

(iii) In the conformal case with $n=2$, there are two $\fG$-invariant bivectors on~$T^*M$,
namely the canonical Poisson bivector and the Poisson bivector
\begin{equation}
\Lambda =
\half
\rg^{ij}\xi_i\xi_j\,\,
\sigma_{k\ell}\,
\frac{\partial}{\partial\xi_k}\wedge\frac{\partial}{\partial\xi_\ell},
\label{extraBivector}
\end{equation}
where $\rg=\rg_{ij}\,dx^idx^j$ represents a conformal class of
(pseudo-)Riemannian metrics and
$\sigma=\half\sigma_{k\ell}\,dx^k\wedge{}dx^\ell$ stands for the surface
element of $(M,\rg)$.
\end{lem}

\begin{proof}
Consider an arbitrary bivector field $W$ on $T^*M$. In any local coordinate system it is of the form
\begin{equation}
W=
W^{ij}(\xi,x)\frac{\partial}{\partial x^i}\wedge\frac{\partial}{\partial{}x^j}+
W_i^j(\xi,x)\frac{\partial}{\partial\xi_i}\wedge\frac{\partial}{\partial{}x^j}+
W_{ij}(\xi,x)\frac{\partial}{\partial\xi_i}\wedge\frac{\partial}{\partial\xi_j},
\label{GeneralBivector}
\end{equation}
where the coefficients
$W^{ij}(\xi,x),W_i^j(\xi,x)$ and $W_{ij}(\xi,x)$ are functions of $x^i,\xi_i$
which are polynomial in~$\xi$.

Choose an adapted coordinate system related to the projective or conformal structure
on $M$ respectively (see Section \ref{ProjSub}).  Since $W$ is $\fG$-invariant, it
commutes with the action of the generators of translations, that is, with the vector
fields $X_i=\partial/\partial x^i$, where $i=1,\ldots,n$. It follows that the
coefficients of $W$ are independent of $x^i$. Furthermore, $W$ is invariant with
respect to the action of the homothety vector field
$X_0=x^i\,\partial/\partial x^i$. The canonical lift of $X_0$ to $T^*M$ is
$L_{X_0}=x^i\,\partial/\partial x^i-\xi_i\,\partial/\partial\xi_i$. One immediately
obtains the following homogeneity conditions:
\begin{enumerate}
\item
the coefficient
$W^{ij}(\xi)$ has to be homogeneous in $\xi$ of degree $-2$,
\item
the coefficient
$W_i^j(\xi)$ has to be homogeneous in $\xi$ of degree $0$,
\item
the coefficient
$W_{ij}(\xi)$ has to be homogeneous in $\xi$ of degree $2$,
\end{enumerate}
so that $W^{ij}(\xi)=0$, while $W_i^j(\xi)$ are constant, and
$W_{ij}(\xi)=W_{ij}^{k\ell}\,\xi_k\xi_\ell$ are quadratic polynomials.
A $\fG$-invariant bivector (\ref{GeneralBivector}) is, therefore, a sum of
two independent $\fG$-invariant bivectors
$W_0=W_i^j\,\partial/\partial\xi_i\wedge\partial/\partial{}x^j$ and
$W_2=W_{ij}^{k\ell}\,\xi_k\xi_\ell\,\partial/\partial\xi_i\wedge\partial/\partial\xi_j$.

Considering now invariance with respect to the linear subgroup of $\fG$ entails that $W_0$ represents an invariant in
$(\bbR^n)^*\otimes\bbR^n$ and $W_2$ an invariant in $\Lambda^2(\bbR^n)^*\otimes{}S^2\bbR^n$ 
with respect to the
standard linear action of $\SL(n,\bbR)$ in the projective case, and~$\SO_0(p,q)$ in the
conformal case. A classical result of invariant theory (see
\cite{Wey,GW}) yields $W_0=c_0\,\Pi$ with $c_0\in\bbC$ and $W_2=0$, except for $n=2$, in the
conformal case, where $W_2=c_2\,\Lambda$ with $c_2\in\bbC$.
Hence, we have proved that the bivectors (\ref{canonicalBivector}),
and (\ref{extraBivector}) for $n=2$ in the conformal
case, are the only bivectors invariant with respect to the affine subgroup of
$\fG$.

To complete the proof, one checks that the bivectors (\ref{canonicalBivector}) and (\ref{extraBivector}) are
invariant with respect to inversions, i.e., the quadratic vector fields in (\ref{ProjVF}) and
(\ref{confGenerators}).
\end{proof}

\begin{lem}
Every $\fG$-invariant Hochschild
$2$-co\-boundary $C_0$ on the
associative commutative algebra $\cS(M)$ is of the form
$C_0=\d{A}$ where $A$ is a $\fG$-invariant linear map on $\cS(M)$.
\label{DifficultLem}
\end{lem}
\begin{proof}
The $2$-coboundary $C_0$ is local thanks to Lemma \ref{LocalityLemma}. This clearly
implies that any $1$-cochain $A$ such that $C_0=\d{A}$ is local, cf.~Proposition
\ref{LemHoch}.

Given a $\fG$-invariant Hochschild
$2$-co\-boundary $C_0=\d{A}$, we will prove that there exists a linear map $\widetilde{A}$
such that
$\d\widetilde{A}=\d{A}$ and
$L_X(\widetilde{A})=0$ for all $X\in\fg$.
Clearly, $\fG$-invariance of $C_0=\d{A}$ implies $L_X(\d{A})=0$ for any
$X\in\fg$. Thus, $\d(L_X({A}))=0$ which means that $L_X({A})$ is a vector field.

A local operator $A$ is a locally given, according to the Peetre theorem \cite{Pee}, by a
differential operator; in an arbitrary coordinate system,
\begin{equation}
A=A^{(0)}+A^{(1)}+A^{(2)}+\cdots+A^{(m)}
\label{A}
\end{equation}
where
\begin{equation}
A^{(i)}=
\sum_{i_1+i_2=i}
A^{s_1\ldots{}s_{i_1}}_{t_1\cdots{}t_{i_2}}(x,\xi)\,
\frac{\partial}{\partial{x^{s_1}}}\cdots\frac{\partial}{\partial{x^{s_{i_1}}}}\;
\frac{\partial}{\partial{\xi_{t_1}}}\cdots\frac{\partial}{\partial{\xi_{t_{i_2}}}}.
\label{Ai}
\end{equation}
Choose a coordinate system adapted to either the projective or the conformal structure.
Consider first the action of the affine Lie subalgebra, $\fh\subset\fg$, that is,
$\fh=\aff(n,\bbR)$ in the projective case and $\fh=\ce(p,q)$ in the conformal case,
introduced in Section \ref{ProjSub}.

For each component~$A^{(i)}$, except for
$A^{(1)}$, one has
$L_X(A^{(i)})=0$, where $X\in\fh$,  since this operator is of the
form~(\ref{Ai}) and thus cannot be a vector field. Put
$\widetilde{A}=A-A^{(1)}$; this operator satisfies $L_X(\widetilde{A})=0$ for all
$X\in\fh$ and, obviously, $\d{\widetilde{A}}=\d{A}$.
In particular, invariance with respect to translations guarantees that
the coefficients in (\ref{Ai}) are independent of~$x$.

In the projective case, an affinely invariant operator $\widetilde{A}$ is of the form
(\ref{AffInvOp}); for the generators
$X_i$ of inversions, $L_{X_i}(\widetilde{A})$ is given by (\ref{LieDerAffInvOp}).
This is a vector field if and only if $P_s=0$ for all $s\geq1$. Hence
$\widetilde{A}=P_0(\cE)$ and thus $L_{X_i}(\widetilde{A})=0$, see
Proposition~\ref{Projinvariants}.

In the conformal case, let us 
rewrite the expression of $\widetilde{A}$ in a different
form, namely
$$
\widetilde{A}=
\widetilde{A}_{(0)}+\widetilde{A}_{(1)}+
\widetilde{A}_{(2)}+\cdots+\widetilde{A}_{(t)}
$$
where $t\leq{}m$ and
$$
\widetilde{A}_{(j)}=
\widetilde{A}^{s_1\ldots{}s_j}\,
\frac{\partial}{\partial{x^{s_1}}}\cdots\frac{\partial}{\partial{x^{s_{j}}}}\,;
$$
each $\widetilde{A}^{s_1\ldots{}s_j}$ is a differential
operator in $\xi$ with polynomial coefficients in $\xi$.
Each term~$\widetilde{A}_{(j)}$ is invariant with
respect to translations and homogeneous in $x$ of degree~$-j$. Invariance 
with respect
to homotheties implies that $\widetilde{A}_{(j)}$ is homogeneous in $\xi$ of
degree~$-j$, that is,
$$
[\cE,\widetilde{A}_{(j)}]=-j\,\widetilde{A}_{(j)}\,.
$$

Let $X_i$ be the $i$-th generator of inversions.
The operator $L_{X_i}(\widetilde{A}_{(j)})$ is homogeneous in~$\xi$ of degree $-j$,
since $L_{X_i}(\cE)=0$, cf.
Proposition~\ref{Confinvariants}. Hence, $L_{X_i}(\widetilde{A})$ is a vector
field only if
$L_{X_i}(\widetilde{A}_{(j)})=0$ for $j\geq2$ since it is polynomial in $\xi$.

Because of its $\fh$-invariance, $\widetilde{A}$ belongs to the ring
generated by the operators
$$
\cE,\quad
\rR_0=\rR\circ\rT,\quad
\rD,\quad
\rG_0=\rG\circ\rT,\quad
\Delta_0=\Delta\circ\rT,
$$
where these operators have been defined in (\ref{sl2}),
(\ref{rh1}) and~(\ref{HomogOpers}). The term $\widetilde{A}_{(1)}$ is then necessarily of the form
$$
\widetilde{A}_{(1)}=\a\,\rD + \b\,\rG_0
$$
where $\a$ and $\b$ are polynomials in $\cE$ and $\rR_0$. A direct computation yields
$$
L_{X_i}(\widetilde{A}_{(1)})=
2\a\left(
\xi_i\rT-2\cE\frac{\partial}{\partial{\xi_i}}-n\frac{\partial}{\partial{\xi_i}}
\right)
- 2\b\left(\rR_0\frac{\partial}{\partial{\xi_i}}+2\xi_i\rT\right).
$$
Every term in this expression, except for $-2n\a\,\partial/\partial{\xi_i}$, is a
differential operator of order~$>1$ for any $\a$ and $\b$. Thus, the right hand
side can be a non-zero vector field only if $\a$ is a non-zero constant. On the other
hand, $-2\b\,\rR_0\,\partial/\partial{\xi_i}$ is, at least, a third-order
term unless $\b$ is zero. But, the remaining terms $2\a\,\xi_i\rT$ and
$-4\a\cE\,\partial/\partial{\xi_i}$ are of order 2 and linearly independent. One
concludes that $\a=0$ and thus $L_{X_i}(\widetilde{A}_{(1)})=0$.

Finally, the term $\widetilde{A}_{(0)}$ is obviously a polynomial in $\cE$ and
$\rR_0$ and, hence,
$L_{X_i}(\widetilde{A}_{(0)})=0$.

We have thus proved that $L_{X_i}(\widetilde{A})=0$ for all $i=1,\ldots,n$.
Lemma \ref{DifficultLem} is proved.
\end{proof}

Let us resort to Lemmas \ref{SkewLem} and
\ref{DifficultLem} to complete the proof. The $\fG$-invariant Hochschild
$2$-cocycle
$C=B_r-B'_r$ is a sum $C=C_1+C_0$.

The symmetric part~$C_0$ is a Hochschild coboundary
and, by Lemma
\ref{DifficultLem} is of the form $C_0=\d{A}$ where $A$ is a $\fG$-invariant
$1$-cochain. This term can be removed by a $\fG$-equivalence map
$\Phi=\Id+(i\hbar)^rA$.

Under the hypotheses of parts (i) and (ii) of Lemma \ref{SkewLem}, the skew-symmetric
part~$C_1$ is proportional to the canonical Poisson bivector, that is, to the first-order
term $B_1$. It can be removed by a reparametrization
$i\hbar\mapsto{}i\hbar+c\,(i\hbar)^r$ for some~$c\in\bbR$.

Theorem \ref{MainThm2} is proved for the first two options,
(i) and (ii), of Lemma \ref{SkewLem}.

In the conformal case and for $n=2$ (part (iii) of Lemma \ref{SkewLem}), the
skew-symmetric part~$C_1$ is a linear combination of the canonical Poisson bivector
$\Pi$ and of the bivector $\Lambda$ in (\ref{extraBivector}). By a reparametrization
map we can remove the canonical Poisson bivector but not the bivector $\Lambda$.

Let us, indeed, show that, if $B_r-B'_r=C_1=k\Lambda$, then necessarily $k=0$.
We associate to the star-products $\star$ and $\star'$ the corresponding
star-commutators
\begin{equation}
[F,G]_\star=\frac{1}{i\hbar}\left(F\star{}G-G\star{}F\right).
\label{starcommutator}
\end{equation}
Since the two star-products are associative, the corresponding star-commutators
satisfy the Jacobi identity. Put
$J_\star(F,G,H)=[F,[G,H]_\star]_\star+\hbox{(cyclic)}$ and consider the difference
$J_\star(F,G,H)-J_{\star'}(F,G,H)$. By assumption, this expression has to be
identically zero. Since the two star-products coincide up to order $r-1$ in $i\hbar$,
this difference is trivially zero up to order $r-2$. Straightforward computation
shows that the $(r-1)$-th order term in the above difference is equal to
$2k[\Pi,\Lambda](F,G,H)$, where $[\Pi,\Lambda]$ is the Schouten bracket of $\Pi$ and
$\Lambda$. Jacobi identities for $\star$ and $\star'$-commutators then lead to
$k\,[\Pi,\Lambda]=0$.

\begin{lem}
\label{Incompatible}
The two Poisson bivectors $\Pi$ and $\Lambda$ are not compatible.
\end{lem}

\goodbreak

\begin{proof}
The Schouten bracket is
\begin{eqnarray*}
[\Pi,\Lambda]
&=&
2\frac{\partial}{\partial\xi_1}\wedge\frac{\partial}{\partial\xi_2}
\wedge\left(
\xi_1\frac{\partial}{\partial{}x^1}+\xi_2\frac{\partial}{\partial{}x^2}
\right)\\[6pt]
&=&
2\Lambda\wedge\frac{\rG}{\rR}
\end{eqnarray*}
where $\rG$ and $\rR$ are as in (\ref{rh1}) and (\ref{sl2}). This expression does not
vanish.
\end{proof}

Thus, the constant $k$ in the above formula has
to vanish. This completes the proof of Part~(iii).

\nobreak
Theorem \ref{MainThm2} is proved.
\end{proof}

\goodbreak

Lemmas \ref{SkewLem}--\ref{Incompatible} can be summarized as the following
\begin{pro}
The second $\fG$-invariant Hochschild cohomology space is
$$
\mathrm{HH}^2_\fG(\cS(M);\cS(M))=
\left\{
\begin{array}{ll}
\bbR^2, & \hbox{in the conformal case for}\ n=2\\
\bbR, & \hbox{otherwise}
\end{array}
\right.
$$
and the cup product in the first instance is non-zero.
\end{pro}

This result could have been derived from Kontsevich's \cite{Kon} or
Fedosov's~\cite{Fed4} classification of equivalence classes of deformations.

\begin{rmk}
{\rm
Theorem \ref{MainThm2} does not guarantee uniqueness of a star-product but of a
class of $\fG$-invariant star-products. Together with Propositions
\ref{Projinvariants} and \ref{Confinvariants} this leads to an explicit
description of all $\fG$-invariant star-products. Indeed, they are all obtained
from the $\fG$-canonical homogeneous star-product by the
equivalence~(\ref{EqUsl}) and reparametrization~(\ref{EqRepar}); the equivalence
map $\Phi$ is given in terms of the $\fG$-invariant operators $\cE$ in the
projective case and $\cE$ and $\rR_0$ in the conformal case. }
\end{rmk}

\subsection{Uniqueness up to $\fG$-equivalence and reparametrization,
$\fG$-covariance and homogeneity}\label{CompareGutt}

In this section we compare our uniqueness theorems with
those obtained for the Moyal star-product in \cite{Gut2}.
The Moyal star-product is the unique, up to
re\-para\-metrization, $(\Sp(2n,\bbR)\ltimes\bbR^{2n})$-invariant
star-product on~$\bbR^{2n}$. It was also proved that it is uniquely selected
within its repa\-rame\-trization class by furthermore requiring
its covariance. The $(\Sp(2n,\bbR)\ltimes\bbR^{2n})$-equivalence class of the
Moyal star-product has a single element since the
$(\Sp(2n,\bbR)\ltimes\bbR^{2n})$-commutant in $\End(C^\infty(\bbR^{2n}))$ is
trivial so that there are no non-zero
invariant Hochschild $2$-coboundaries.

One may wonder if in our present setting $\fG$-covariance
plays a similar role, namely, that of an extra condition
that selects the canonical $\fG$-invariant star-product
of Section~\ref{HomogUniq} within its reparametrization and
$\fG$-equivalence classes described in
Section~\ref{UniqueUpGequiReparam}. The answer is negative; however we have

\begin{pro}\label{GCovUnique}
If two $\fG$-invariant and $\fG$-covariant star-products on $\cS(M)$ are
 equivalent up to reparametrization, then they coincide.
\end{pro}
\begin{proof}
Let $\star$ and $\star'$ be two $\fG$-invariant and $\fG$-covariant
star-products on $\cS(M)$ belonging to the same reparametrization
class. Their $\fG$-covariance translates into
(see (\ref{covStarProdDef}))~:
\begin{equation}\label{GCov2}
J_X\star{}J_Y-J_Y\star{}J_X=i\hbar\,\{J_X,J_Y\}=
J_X\star'{}J_Y-J_Y\star'{}J_X
\end{equation}
for all $X,Y\in\fg$. On the other hand reparametrization
equivalence means that there exist a formal power series (\ref{EqRepar})
such that
$$
F\star'{}G=\sum_{r\geq0}{(i\hbar)^r{}B_r'(F,G)}=
\sum_{r\geq0}\left(\mu(i\hbar)\right)^r{}B_r(F,G).
$$
Using this equation, one rewrites the right hand side of (\ref{GCov2})
in terms of $\star$, with $\mu(i\hbar)$ as deformation
parameter. Now, using the left hand side of (\ref{GCov2}) one gets
$\mu(i\hbar)=i\hbar$, from which the conclusion follows.
\end{proof}

An analog of the above statement, where the reparametrization
equivalence is replaced by $\fG$-equivalence, does not hold.
Indeed, one shows using an argument similar to the one in the
above proof, that two $\fG$-invariant and $\fG$-covariant star-products on
$\cS(M)$ in the same $\fG$-equivalence class, do not
necessarily coincide. So, covariance does not play the same
role for $\fG$ as it does for $\Sp(2n,\bbR)\ltimes\bbR^{2n}$. However,
a simple verification shows that, for the Moyal star-product, homogeneity 
has exactly the same effect as $(\Sp(2n,\bbR)\ltimes\bbR^{2n})$-covariance.
Hence, the $\fG$-canonical and the Moyal star-products are uniquely
determined by two simple conditions, namely, invariance and homogeneity.

\section{Explicit formula for the projectively-invariant star-product}
\label{ExplicitProjSection}

In this section we compute the explicit formula of the canonical homogeneous
projectively-invariant star-product. This solves a problem raised in \cite{AB}.

Projective invariance will be dealt with in two stages. We first consider
invariance with respect to an affine subgroup $\Aff(n,\bbR)$ of $\SL(n+1,\bbR)$ and
determine the affine-invariant bilinear operators on $\cS(\bbR\rP^n)$. Those will
be used to write down an Ansatz for the star-product we are
looking for.  We will then enforce full projective invariance by further demanding
that inversions preserve the star-product. This will give rise to the equations
(\ref{InversionsInvariance1}) and (\ref{InversionsInvariance2}) below. Another system of
equations will arise from the associativity requirement (see
(\ref{AssociativityCondition})).
The unique solution of the complete system of equations will be given explicitly at the
end of this section.

\subsection{Autonomous derivation from the invariance principle}

We need to classify the bilinear $\Aff(n,\bbR)$-invariant differential
operators on $\cS(\bbR^n)$. For that purpose, let us resort to the natural
isomorphism
\begin{equation}
\label{IdentBin}
\cS(\bbR^n)\otimes\cS(\bbR^n)
\cong
\cS(\bbR^{2n})
\end{equation}
and denote by $(x,\xi,y,\eta)$ the natural coordinate system on
$T^*\bbR^{n}\times{}T^*\bbR^{n}$. The operators
of divergence with respect to the first and the second
arguments
\begin{equation}
\rD_{x\xi}(F,G)=D(F)\,G,
\qquad
\qquad
\rD_{y\eta}(F,G)=F\,D(G),
\label{Div}
\end{equation}
where $D$ is as in (\ref{ED}), and the operators of contraction
\begin{eqnarray}
\rD_{x\eta}(F,G)
&=&
{\frac{\partial}{\partial{}x^i}\frac{\partial}{\partial{}\eta_i}
F(\xi,x)G(\eta,y)
}
\Big\vert_{\eta=\xi,y=x},\\[8pt]
\rD_{y\xi}(F,G)
&=&
{\frac{\partial}{\partial{}y^i}\frac{\partial}{\partial{}\xi_i}
F(\xi,x)G(\eta,y)
}
\Big\vert_{\eta=\xi,y=x}
\label{Dcontract}
\end{eqnarray}
are obviously $\Aff(n,\bbR)$-invariant differential operators.
Restricting ourselves to homogeneous components, we get the following
\begin{pro}
\label{AffPro}
Every bilinear differential operator
\begin{equation}
\label{Bilin1}
\cS_k(\bbR^n)\otimes\cS_\ell(\bbR^n)\to\cS_m(\bbR^n)
\end{equation}
invariant with respect to the action of $\Aff(n,\bbR)$, is a homogeneous
polynomial in $\rD_{\xi{}x},\rD_{\xi{}y},\rD_{\eta{}x}$ and $\rD_{\eta{}y}$
of degree $k+\ell-m$.
\end{pro}

This enables us to write the most general $\Aff(n,\bbR)$-invariant bilinear
operation $\cS(\bbR^n)\otimes\cS(\bbR^n)\to\cS(\bbR^n)[\hbar]$.  According to
Theorem \ref{MainThm2}, we will express it as a termwise homogeneous formal
series which, when restricted to
$\cS_k(\bbR^n)\otimes\cS_\ell(\bbR^n)$, takes the form
\begin{equation}
F\star{}G=\sum_{r=0}^\infty{(i\hbar)^rB_r^{k,\ell}(F,G)}
\label{Ansatz1}
\end{equation}
where $B_r^{k,\ell}$ is a bidifferential operator, homogeneous of degree $r$ in
$\rD_{\xi{}x},\rD_{\xi{}y},\rD_{\eta{}x},\rD_{\eta{}y}$, viz
\begin{equation}
B_r^{k,\ell}(F,G)(\xi,x)
=
\sum_{\a+\b+\c+\d=r}{\!\!
B_{\a,\b,\c,\d}^{k,\ell}\,
\rD_{\xi{}y}^\a\,\rD_{\eta{}x}^\b\,\rD_{\xi{}x}^\c\,\rD_{\eta{}y}^\d\,F(\xi,x)G(\eta,y)
}
\Big\vert_{\eta=\xi,y=x}
\label{Ansatz2}
\end{equation}
with constant coefficients $B_{\a,\b,\c,\d}^{k,\ell}$.

\goodbreak

Since we seek a star-product, we have to impose
\begin{equation}
B_{0,0,0,0}^{k,\ell}=1
\label{B0}
\end{equation}
and
\begin{equation}
B_{1,0,0,0}^{k,\ell}=-B_{0,1,0,0}^{k,\ell}=\half
\label{B1}
\end{equation}
in order to get the multiplication and Poisson bracket as the first
two terms as in equation (\ref{starproduct}).

Expressions (\ref{Ansatz1}) and (\ref{Ansatz2}) constitute our Ansatz
for an
$\SL(n+1,\bbR)$-invariant star-product on $T^*\bbR\rP^n$. It now remains to 
impose to
the operation (\ref{Ansatz1}) the following conditions: (i) invariance with 
respect to
inversions, and (ii) associativity.

\subsection{Projective invariance}

Let $X_i=x^ix^j\partial_{x^j}$ be the $i$-th generator of inversions. Denote by
$$
L_{X_i}
=
x^ix^j\frac{\partial}{\partial{x^j}}
-
x^j\xi_j\frac{\partial}{\partial{\xi_i}}
-
x^i\xi_j\frac{\partial}{\partial{\xi_j}}
+
y^iy^j\frac{\partial}{\partial{y^j}}
-
y^j\eta_j\frac{\partial}{\partial{\eta_i}}
-
y^i\eta_j\frac{\partial}{\partial{\eta_j}}
$$
its canonical lift to $T^*(\bbR^{2n})$.

Invariance with respect to inversions translates into the following equations
\begin{equation}
\sum_{\a+\b+\c+\d=r}\!\!
B_{\a,\b,\c,\d}^{k,\ell}\,
\left[L_{X_i},
\rD_{\xi{}y}^\a\,\rD_{\eta{}x}^\b\,\rD_{\xi{}x}^\c\,\rD_{\eta{}y}^\d
\right]\Big\vert_{\eta=\xi,y=x}
=0
\label{InversionsCommRels}
\end{equation}
at each order $r\in\bbN$. The latter yield the
following system of equations
\begin{eqnarray}\label{InversionsInvariance1}
\lefteqn{(\a+1)(\a+\d-\ell)B^{k,\ell}_{\a+1,\b,\c,\d}
+(\b+1)(\b+\d-\ell)B^{k,\ell}_{\a,\b+1,\c,\d}=}\nonumber\\[12pt]
&&(\c+1)(n+2k-\c-1)B^{k,\ell}_{\a,\b,\c+1,\d}
+(\a+1)(\b+1)B^{k,\ell}_{\a+1,\b+1,\c,\d-1}
\end{eqnarray}
and
\begin{eqnarray}\label{InversionsInvariance2}
\lefteqn{(\b+1)(\b+\c-k)B^{k,\ell}_{\a,\b+1,\c,\d}
+(\a+1)(\a+\c-k)B^{k,\ell}_{\a+1,\b,\c,\d}=}\nonumber\\[12pt]
&&(\d+1)(n+2\ell-\d-1)B^{k,\ell}_{\a,\b,\c,\d+1}
+(\a+1)(\b+1)B^{k,\ell}_{\a+1,\b+1,\c-1,\d}.
\end{eqnarray}

\subsection{Associativity}

If $F\in\cS_k(\bbR^n)$, $G\in\cS_\ell(\bbR^n)$, and $H\in\cS_m(\bbR^n)$ the
associativity condition takes the form
\begin{equation}
\sum_{j=0}^r{
B_{r-j}^{k,\ell+m-j}(F,B_j^{\ell,m}(G,H))
}
=
\sum_{j=0}^r{
B_{r-j}^{k+\ell-j,m}(B_j^{k,\ell}(F,G),H)
}
\label{Associativity1}
\end{equation}
for all $r\in\bbN$.
Equation (\ref{Associativity1}) then reads
\begin{equation}\label{AssociativityCondition}
\begin{array}{rcl}
&\displaystyle\sum_{j=0}^r
&\displaystyle\sum_{\a+\b+\c+\d=r-j}
B^{k,\ell+m-j}_{\a,\b,\c,\d}
(\rD_{\xi{}y}+\rD_{\xi{}z})^\a
(\rD_{\eta{}x}+\rD_{\zeta{}x})^\b\times\\[12pt]
&&
\rD_{\xi{}x}^\c
(\rD_{\eta{}y}+\rD_{\eta{}z}+\rD_{\zeta{}y}+\rD_{\zeta{}z})^\d{}\times\nonumber\\[12pt]
&&
\displaystyle\sum_{\a'+\b'+\c'+\d'=j}
B^{\ell,m}_{\a',\b',\c',\d'}
\rD_{\eta{}z}^{\a'}\rD_{\zeta{}y}^{\b'}\rD_{\eta{}y}^{\c'}\rD_{\zeta{}z}^{\d'}=
\nonumber\\[12pt]
&\displaystyle\sum_{j=0}^r
&\displaystyle\sum_{\a+\b+\c+\d=r-j}
B^{k+\ell-j,m}_{\a,\b,\c,\d}
(\rD_{\xi{}z}+\rD_{\eta{}z})^\a
(\rD_{\zeta{}x}+\rD_{\zeta{}y})^\b\times\nonumber\\[12pt]
&&
(\rD_{\xi{}x}+\rD_{\xi{}y}+\rD_{\eta{}x}+\rD_{\eta{}y})^\c
\rD_{\zeta{}z}^\d{}\times\nonumber\\[10pt]
&&
\displaystyle\sum_{\a'+\b'+\c'+\d'=j}
B^{k,\ell}_{\a',\b',\c',\d'}
\rD_{\xi{}y}^{\a'}\rD_{\eta{}x}^{\b'}\rD_{\xi{}x}^{\c'}\rD_{\eta{}y}^{\d'}.
\end{array}
\end{equation}

\subsection{Explicit solution of the system}\label{ExplicitStar}

We solve the system of equations (\ref{InversionsInvariance1}),
(\ref{InversionsInvariance2}) and (\ref{AssociativityCondition}),
by first determining the
components $B_{\a,\b,0,0}^{k,\ell}$, then $B_{\a,\b,\c,0}^{k,\ell}$ and,
finally, the full expression  $B_{\a,\b,\c,\d}^{k,\ell}$.

\subsubsection{First stage}

Identifying in the associativity equation (\ref{AssociativityCondition})
the coefficients of the monomials
$\rD_{\xi{}z}^{r-j}\rD_{\zeta{}x}^{j}$, one readily finds
$B_{r-j,j,0,0}^{k,\ell+m}=B_{r-j,j,0,0}^{k+\ell,m}$. Thus,
$B_{\a,\b,0,0}^{k,\ell}$ depends only on
$k+\ell$; we write
\begin{equation}
B_{\a,\b,0,0}^{k,\ell}=C_{\a,\b}(k+\ell).
\label{Bab00}
\end{equation}

Using again (\ref{AssociativityCondition}), we identify the coefficients of the
monomials
$\rD_{\xi{}z}^{r-j-1}\rD_{\zeta{}x}^{j}\rD_{\eta{}x}$ and
$\rD_{\xi{}z}^{r-j-1}\rD_{\zeta{}x}^{j}\rD_{\xi{}y}$, respectively, to get the
following system
\begin{eqnarray*}
(j+1)B_{r-j-1,j+1,0,0}^{k,\ell+m}
&=&
B_{r-j-1,j,1,0}^{k+\ell,m}+B_{r-j-1,j,0,0}^{k+\ell,m}B_{0,1,0,0}^{k,\ell},\\[12pt]
(r-j)B_{r-j,j,0,0}^{k,\ell+m}
&=&
B_{r-j-1,j,1,0}^{k+\ell,m}+B_{r-j-1,j,0,0}^{k+\ell-1,m}B_{1,0,0,0}^{k,\ell}.
\end{eqnarray*}
Resorting to the invariance equation (\ref{InversionsInvariance1}) for $\a=r-j-1$,
$\b=j$, and $\c=\d=0$, we obtain the supplementary equation
$$
(r-j)(r-j-\ell-1)B^{k,\ell}_{r-j,j,0,0}
+(j+1)(j-\ell)B^{k,\ell}_{r-j-1,j+1,0,0}
-(n+2k-1)B^{k,\ell}_{r-j-1,j,1,0}=0.
$$
The previous three equations together with (\ref{B1}) and (\ref{Bab00}) imply
$$
(r-j)(r-n-2k)C_{r-j,j}(k)+\half(n+2k-2j-1)C_{r-j-1,j}(k-1)=0.
$$
The latter equation, supplemented with (\ref{B0}), yields then
\begin{equation}
B_{\a,\b,0,0}^{k,\ell}
=
\frac{(-1)^\b}{(\a+\b)!}
\frac{{\half(n-1)+k+\ell-\b\choose\a}{\half(n-1)+k+\ell-\a\choose\b}}%
{{n+2k+2\ell-\a-\b\choose\a+\b}}.
\label{SolBab00}
\end{equation}

\subsubsection{Second stage}

Here we only use the first invariance equation (\ref{InversionsInvariance1}) with $\d=0$.
Long but straightforward calculations lead to
\begin{equation}
B_{\a,\b,\c,0}^{k,\ell}
=
\frac{1}{\c!(n+2k-\c)_\c}
\sum_{r+s=\c}{
{\c\choose{}r}(\a+1)_r(\b+1)_s(\a-\ell)_r(\b-\ell)_s
B_{\a+r,\b+s,0,0}^{k,\ell}
}
\label{SolBabc0}
\end{equation}
where the last term is as in (\ref{SolBab00}).

\subsubsection{Last stage}

A reverse iterative computation on $\d$ using the second invariance equation
(\ref{InversionsInvariance2}) finally leads to the sought for result
\begin{eqnarray}
B_{\a,\b,\c,\d}^{k,\ell}
&=&
\frac{1}{\d!(n+2\ell-\d)_\d}
\sum_{r+s+t=\d}
(-1)^s
{\d\choose{}{r,s,t}}\times\nonumber\\[12pt]
& &
(\a+1)_r(\a+1)_s(\b+1)_s(\b+1)_t(\a+\c-k)_r(\b+\c-k)_t\times\nonumber\\[12pt]
& &
B_{\a+r+s,\b+s+t,\c-s,0}^{k,\ell}
\label{SolBabcd}
\end{eqnarray}
where the first line contains the trinomial coefficient and the last one is
given by~(\ref{SolBabc0}).

\subsection{Symmetry condition}

\begin{pro}
The symmetry condition C2
translates for the Ansatz (\ref{Ansatz1})--(\ref{Ansatz2}) into
\begin{equation}
B_{\a,\b,\c,\d}^{k,\ell}=(-1)^{\a+\b+\c+\d}\,B_{\b,\a,\d,\c}^{\ell,k}.
\label{SymmCond2}
\end{equation}
\end{pro}
\begin{proof}
If $F\in\cS_k(\bbR^n)$, and $G\in\cS_\ell(\bbR^n)$, we immediately get from
Condition C2 that
$$
B_r^{k,\ell}(F,G)=(-1)^r\,B_r^{\ell,k}(G,F).
$$
Then, a change of dummy variables in (\ref{Ansatz2}) completes the proof.
\end{proof}

It turns out that our star-product given by (\ref{Ansatz1}), (\ref{Ansatz2}) and
(\ref{SolBabcd}) automatically satisfies the symmetry condition (\ref{SymmCond2}).
Although this is not transparent from the expression (\ref{SolBabcd}), it is however a
direct consequence of Proposition \ref{SymmStarQhalf} and Theorem~\ref{MainThm1}.

\section{Conclusion, discussion and outlook}


In this work we have proved the existence and uniqueness of a canonical
$\fG$-invariant star-product on $T^*M$ for $\fG=\SL(n+1,\bbR)$
(resp. $\fG=\SO_0(p+1,q+1)$ and $M=\bbR{}\rP^n$
(resp.~$(S^p\times S^q)/\bbZ_2$). We have, moreover, given an explicit formula
for the canonical projectively invariant star-product. For both geometries, the
canonical star-product so obtained is symmetric, homogeneous, strongly
$\fG$-invariant (hence $\fG$-covariant), but not differential.
These properties, except for the
last one, are shared with the Moyal star-product on $\bbR^{2n}$.


Theorem \ref{MainThm1} shows that the
homogeneity condition supplementing $\fG$-invariance
uniquely determines the canonical  $\fG$-invariant star-product
on $\cS(M)$.  Likewise,  the Moyal star-product
is also uniquely specified by $(\Sp(2n,\bbR)\ltimes\bbR^{2n})$-invariance and
homogeneity.  This allows us to draw a
parallel between our canonical $\fG$-invariant star-product and
Moyal's, namely, they are uniquely determined by the same two
simple conditions~: invariance and homogeneity. Of course, this parallel
is far from complete, since, for instance, $\fG$ and
$\Sp(2n,\bbR)\ltimes\bbR^{2n}$ do not have the same geometric status; the
action of the former on $T^*M$ is lifted from that on $M$, which is not the case
for the latter.

Furthermore, it is clear that, for the projective and the conformal cases, there
is no $\fG$-invariant (symplectic) connection on $T^*M$,
since~$\fG$ does not act on the bundle of linear frames of $T^*M$.
Hence, no Fedosov \cite{Fed3} canonical $\fG$-invariant star-product
can be constructed. Besides, Fedosov's construction would have led to
a star-product given by bidifferential operators.

The generalization of the existence and uniqueness theorems
for projectively/con\-fo\-rmally invariant star-products on $T^*M$ in the case of a
non-flat projective/con\-fo\-rmal connection on $M$ remains an open problem.
In a recent work \cite{Bor}, Bordemann has taken a significant step
in this new direction, by investigating the projectively equivariant
quantization on a cotangent bundle of a manifold with a non flat
projective structure (see also \cite{DO1} and \cite{Bo}).  Note also, that since the canonical
star-products  studied in this work may be considered as the projective/conformal analogs of 
the Moyal star-product, they may play a similar role as the latter in a 
construction \`a la Fedosov of a star-product on a symplectic manifold with a
Cartan projective/conformal symplectic connection.

In the case $n\geq2$, let us mention that the explicit form of the
conformally invariant star-product is, so far, out of reach. This was
already the situation for the conformally equivariant quantization map
\cite{DLO} (see also \cite{DO1}).

In the conformal case with $n=2$, Theorem \ref{MainThm2} holds for star-products of
the form (\ref{starproduct}) with the standard Poisson bracket on $T^*M$ as
first-order term. However, one could easily construct, in this case, another
$\fG$-invariant star-product with the Poisson bracket (\ref{extraBivector}) as
first-order term. It would be interesting to give a physical status to this
second, somewhat ``exotic'', star-product.

In the case of dimension $n=1$, our results are related to earlier work by
Cohen, Manin and Zagier~\cite{CMZ}. The projective and the conformal algebras
are, in this case, isomorphic to $\Sl(2,\bbR)$. Moreover, the canonical
projectively and conformally invariant star-products coincide by uniqueness  and
thus the explicit formul\ae\  given in Section~\ref{ExplicitStar} correspond to the
one obtained in
\cite{CMZ} for~$\l=\half$.


\end{document}